\documentclass[12pt,letterpaper,final,twoside,leqno]{amsart}
\usepackage{amsmath,amssymb,amsthm,amsfonts,amscd,amsopn}
\usepackage[T1]{fontenc}
\usepackage{times}
\usepackage{epsfig}
\usepackage{xspace}
\usepackage[all]{xy}
\xyoption{dvips}
\newcommand\mypagesize{
\textwidth= 6.25in
\textheight=8.75in
\voffset-.5in
\hoffset-.7in
}
\makeatletter
\mypagesize
\DeclareFontFamily{OMS}{rsfs}{\skewchar\font'60}
\DeclareFontShape{OMS}{rsfs}{m}{n}{<-5>rsfs5 <5-7>rsfs7 <7->rsfs10 }{}
\DeclareSymbolFont{rsfs}{OMS}{rsfs}{m}{n}
\DeclareSymbolFontAlphabet{\scr}{rsfs}

\newcommand{\sA}{\scr{A}}

\newcommand{\sE}{\scr{E}}
\newcommand{\sF}{\scr{F}}

\newcommand{\sH}{\scr{H}}

\newcommand{\sL}{\scr{L}}

\newcommand{\sO}{\scr{O}}
\newcommand{\sP}{\scr{P}}

\newcommand{\bC}{\mathbb{C}}

\newcommand{\bP}{\mathbb{P}}

\newcommand{\cP}{\mathcal{P}}

\newcommand{\cW}{\mathcal{W}}

\newcommand{\cZ}{\mathcal{Z}}

\newcommand{\al}{\alpha}
\newcommand{\be}{\beta}
\newcommand{\ga}{\gamma}

\newcommand{\into}{\hookrightarrow}

\newcommand{\leteq}{\,\colon\!\!=}

\renewcommand{\l}{\scr L}

\newcommand{\N}{\mathbb {N}}
\newcommand{\PP}{\mathbb {P}}

\newcommand{\tensor}{\otimes}

\newcommand\hypush{R\phi_{\cdot *}}
\def\rpforward#1{R{#1}_*}

\DeclareMathOperator{\alb}{{alb}}
\DeclareMathOperator{\Alb}{{Alb}}

\DeclareMathOperator{\id}{{id}}
\DeclareMathOperator{\im}{{im}}

\DeclareMathOperator{\ob}{{Ob}}

\DeclareMathOperator{\Pic}{Pic}

\DeclareMathOperator{\rank}{{rank}}

\def\spec#1.#2.{{\mathbf S\mathbf p\mathbf e\mathbf c}_{#1}#2}
\def\proj#1.#2.{{\mathbf P\mathbf r\mathbf o\mathbf j}_{#1}\sum#2}
\def\ring#1.{\scr O_{#1}}
\def\map#1.#2.{#1 \to #2}

\def\longmap#1.#2.{#1 \longrightarrow #2}
\def\factor#1.#2.{\left. \raise 2pt\hbox{$#1$} \right/
\hskip -2pt\raise -2pt\hbox{$#2$}}
\def\pe#1.{\mathbb P(#1)}
\def\pr#1.{\mathbb P^{#1}}

\def\Oxi#1.{{\Omega}^{#1}_{X}}
\def\Ox{\Omega_{X}}
\def\ox{\omega_{X}}
\def\Oxu{\underline{\Omega}_{X}}
\def\Oxph{{\Omega}_{X_\cdot}^p}
\def\Oxh#1{{\Omega}_{X_\cdot}^{#1}}
\def\OxDi#1.{{\Omega}_X^{#1}(\log{D})}
\def\Oxx#1.#2.{{\Omega}_X^{#1}(\log{#2})}

\def\Oxs#1.{{\Omega}_{X/S}^{#1}(\log{D})}

\def\Os#1.{{\Omega}_S^{#1}(\log{\Delta})}
\def\Oss#1.#2.{{\Omega}_S^{#1}(\log{#2})}

\def\oy{\omega_{Y}}
\def\foy{\phi^{^\#}\!\Omega_Y}
\def\Oyu{\underline{\Omega}_{Y}}
\def\Oyp{\underline{\Omega}_{Y}^p}

\def\coh#1.#2.#3.{\dsize H^{#1}(#2,#3)}
\def\dimcoh#1.#2.#3.{\dsize h^{#1}(#2,#3)}
\def\hypcoh#1.#2.#3.{\dsize \mathbb H_{\vphantom{l}}^{#1}(#2,#3)}
\def\loccoh#1.#2.#3.#4.{\dsize H^{#1}_{#2}(#3,#4)}
\def\dimloccoh#1.#2.#3.#4.{\dsize h^{#1}_{#2}(#3,#4)}
\def\lochypcoh#1.#2.#3.#4.{\dsize \mathbb H^{#1}_{#2}(#3,#4)}
\def\ses#1.#2.#3.{0  \longrightarrow  #1   \longrightarrow 
 #2 \longrightarrow #3 \longrightarrow 0} 
\def\sesshort#1.#2.#3.{0
 \rightarrow #1 \rightarrow #2 \rightarrow #3 \rightarrow 0}
\def\dist#1.#2.#3.{  #1   \longrightarrow 
 #2 \longrightarrow #3 \stackrel{+1}{\longrightarrow} } 
\def\CDdist#1.#2.#3.{  #1   @>>>  #2  @>>>   #3 @>+1>> }  
\def\shortses#1.#2.#3.{0  \rightarrow  #1   \rightarrow 
 #2  \rightarrow   #3 \rightarrow  0}
\def\shortdist#1.#2.#3.{  #1   \rightarrow 
 #2  \rightarrow   #3 \stackrel{+1}{\rightarrow} }  
\def\ddist#1.#2.#3.#4.#5.#6.{\CD
#1 @>>> #2 @>>> #3 @>+1>> \\
@VVV @VVV @VVV \\
#4 @>>> #5 @>>> #6 @>+1>> 
\endCD}
\def\ddistun#1.#2.#3.#4.#5.#6.{\CD
#1 @>>> #2 @>>> #3 @>+1>> \\
@. @VVV @VVV  \\
#4 @>>> #5 @>>> #6 @>+1>> 
\endCD}
\def\Iff#1#2#3{
\hfil\hbox{\hsize =#1
\vtop{\noin #2}
\hskip.5cm 
\lower.5\baselineskip\hbox{$\Leftrightarrow$}\hskip.5cm
\vtop{\noin #3}}\hfil\medskip}
\def\myoplus#1.#2.{\underset #1 \to {\overset #2 \to \oplus}}

\def\resto#1{\hbox{\hbox{$\vert{}_{\dsize_{#1}}$}}}
\def\leftresto{\left.}
\def\rightresto#1{\right\vert{}_{{}_{\dsize_{#1}}}}

\def\ww#1.#2.{\curlywedge
_{#1}^{#2}}
\def\wws#1.#2.{\scriptstyle\ww#1.#2.}
\def\nattr{\theta}
\def\fh#1.#2.{F^{#1} R^{#2}}
\def\epq#1.#2.#3.{E^{#1,#2}_{#3}}
\def\dd#1.#2.#3.{d^{\, #1,#2}_{#3}}
\def\kk#1.#2.#3.{K^{\, #1,#2}_{#3}}
\def\ki#1.#2.#3.{I^{\, #1,#2}_{#3}}
\def\ff#1.#2.{\mathfrak F_{#1}^{#2}}
\def\ffb#1.#2.{\boxed{\mathfrak F_{#1}^{#2}}}
\def\ffbx#1.#2.#3.{\boxed{\mathfrak F_{#1}^{#2}(#3)}}
\def\ffs#1.#2.{\scriptstyle\mathfrak F_{#1}^{#2}}
\def\ffsb#1.#2.{\boxed{\ffs#1.#2.}}
\def\ffsbx#1.#2.#3.{\boxed{\scriptstyle\mathfrak F_{#1}^{#2}(#3)}}
\def\fa#1.#2.{\mathfrak f_{#1}^{#2}}
\def\fab#1.#2.{\boxed{\mathfrak f_{#1}^{#2}}}
\def\ffa#1.#2.{\mathfrak F\mathfrak i\mathfrak l\mathfrak t_{#1}^{#2}}
\def\ffab#1.#2.{\boxed{\mathfrak F\mathfrak i\mathfrak l\mathfrak t_{#1}^{#2}}}
\def\rfi#1.#2.{R^{#1}\functor(#2\,)}

\def\al#1.#2.#3.{\alpha_{#1,#2}^{#3}}
\def\be#1.#2.{\beta_{#1,#2}}
\def\ga#1.#2.{\gamma_{#1,#2}}
\def\yy#1.#2.#3.#4.{y_{#1,#2}^{#3}(#4)}
\def\zz#1.#2.#3.#4.{z_{#1,#2}^{#3}(#4)}
\newcounter{lastyear}
\addtocounter{lastyear}{-1}

\newcommand\noin{\noindent}

\newcommand\dsize{\displaystyle}
\newcommand\sh{\subheading}

\newcommand\myref[1]{(\ref{#1})}

\def\vol#1 {{\bf #1}\ } 
\def\yr#1 {\rm (#1)\ } 

\newcommand\te{there exist}
\newcommand\st{such that}

\newtheoremstyle{bozont}{3pt}{3pt}%
     {\itshape}
     {}
     {\bfseries}
     {.}
     {.5em}
     {\thmname{#1}\thmnumber{ #2}\thmnote{ #3}}
\newtheoremstyle{bozont-remark}{3pt}{3pt}%
     {}
     {}
     {\sc}
     {.}
     {.5em}
     {\thmname{#1}\thmnumber{ #2}\thmnote{ #3}}
\newtheoremstyle{bozont-def}{3pt}{3pt}%
     {}
     {}
     {\bfseries}
     {.}
     {.5em}
     {\thmname{#1}\thmnumber{ #2}\thmnote{ #3}}
\newtheoremstyle{bozont-reverse}{3pt}{3pt}%
     {\itshape}
     {}
     {\bfseries}
     {.}
     {.5em}
     {\thmnumber{#2}\thmname{ #1}\thmnote{ #3}}
\newtheoremstyle{bozont-remark-reverse}{3pt}{3pt}%
     {}
     {}
     {\sc}
     {.}
     {.5em}
     {\thmnumber{#2}\thmname{ #1}\thmnote{ #3}}
\newtheoremstyle{bozont-def-reverse}{3pt}{3pt}%
     {}
     {}
     {\bfseries}
     {.}
     {.5em}
     {\thmnumber{#2}\thmname{ #1}\thmnote{ #3}}
\theoremstyle{bozont}    
\newtheorem{proclaim}{Theorem} [section]
\newtheorem{thm}[proclaim]{Theorem}
\newtheorem{cor}[proclaim]{Corollary} 
\newtheorem{lem}[proclaim]{Lemma} 
\newtheorem{prop}[proclaim]{Proposition} 
\newtheorem{conj}[proclaim]{Conjecture}
\newtheorem{proclaim-special}[proclaim]{\specialthmname}
\newenvironment{proclaimspecial}[1]
     {\def\specialthmname{#1}\begin{proclaim-special}}
     {\end{proclaim-special}}
\theoremstyle{bozont-remark}
\newtheorem{rem}[proclaim]{Remark}
\newtheorem{notation}[proclaim]{Notation}

\newtheorem*{SubHeading*}{\SubHeadingName}

\newenvironment{demo}[1]
     {\def\SubHeadingName{#1}\begin{SubHeading}}
     {\end{SubHeading}}
\newenvironment{demo*}[1]
     {\def\SubHeadingName{#1}\begin{SubHeading*}}
     {\end{SubHeading*}}

\theoremstyle{bozont-def}    
\newtheorem{defn}[proclaim]{Definition}

\theoremstyle{bozont-reverse}    
\newtheorem{proclaimr}{Theorem} [section]

\newtheorem{proclaimr-special}[proclaimr]{\specialthmname}

\theoremstyle{bozont-remark-reverse}

\newtheorem{SubHeading}[proclaim]{\SubHeadingName}
\newtheorem{SubHeadingr}[proclaimr]{\SubHeadingName}

\theoremstyle{bozont-def-reverse}

\numberwithin{equation}{proclaim}
\numberwithin{figure}{section} 

\def\Section#1
        { 
        \setcounter{firstsubsection}{0}
        \setcounter{proclaim}{0}
        \section{#1}
        \numberwithin{proclaim}{section}
        \rm\noindent
        }
\newcounter{firstsubsection}
\def\Subsection#1
        {
\ifnum \value{firstsubsection}=0 
        \setcounter{subsection}{\value{proclaim}}
        \setcounter{firstsubsection}{1}
\fi
\noin\subsection{#1} 
        \setcounter{proclaim}{0}
        \numberwithin{proclaim}{subsection}
        \numberwithin{equation}{subsection}
        \numberwithin{figure}{subsection} 
        }

\newcounter{rosternumber}\numberwithin{rosternumber}{proclaim}
\newcounter{rosternumberr}\numberwithin{rosternumberr}{proclaimr}
\newenvironment{roster} {\begin{enumerate}
    
}
  {\end{enumerate}}

\def\rostitem{ \refstepcounter{rosternumber}
        \item[(\therosternumber)]
        }

\newcommand\subheading[1]{\smallskip\noindent{{\bf #1.}\ }}
\makeatother

\begin{document}

\title{Holomorphic one-forms on varieties of general type}
\author{Christopher D.\ Hacon and S\'andor J.\ Kov\'acs} 
\begin{abstract}
  It has been conjectured that varieties of general type do not admit nowhere
  vanishing holomorphic one-forms. We confirm this conjecture for smooth minimal
  varieties and for varieties whose Albanese variety is simple.
\end{abstract}
\thanks{Christopher Hacon and was partially supported by NSA grant MDA904-03-1-0101
  and by a grant from the Sloan Foundation.  \newline\indent S\'andor Kov\'acs was
  partially supported by NSF grant DMS-0092165 and a Sloan Research Fellowship.}

\address{University of Utah, Department of Mathematics, 155 South 1400 East, Room
  233, Salt Lake City, UT 84112, USA} \email{hacon@math.utah.edu}
\urladdr{http://www.math.utah.edu/$\sim$hacon} \address{University of Washington,
  Department of Mathematics, 354350, Seattle, WA 98195, USA}
\email{kovacs@math.washington.edu\xspace}
\urladdr{http://www.math.washington.edu/$\sim$kovacs\xspace}
\maketitle
\def\Oxsui#1.#2.{\underline{\mathfrak Q}^{#1}_{\,\nattr_{#2}}} 
\def\foy{\phi^{^\#}\!\Omega_Y}
\newcommand\ot{{\otimes}}
\nocite{Kovacs97c}

\Section{Introduction}
\noin The impact of zeros of vector fields on the geometry of the underlying variety
has been studied extensively, cf.\ \cite{Bott67}, \cite{Baum-Bott70},
\cite{Baum-Bott72}, \cite{CL73}, \cite{CHK73}, \cite{Car-Lie77}, \cite{ACL81},
\cite{Aky-Car83}, \cite{ACL86}. For instance, it is known that the existence of a
nowhere zero vector field on a compact complex manifold implies that all of its
characteristic numbers vanish.

Carrell asked whether something similar is implied by the existence of a nowhere
vanishing holomorphic one form. He proved that this is the case for surfaces, namely
if $S$ is a compact complex surface admitting a nowhere vanishing holomorphic one
form, then $c_1(S)^2$ and $c_2(S)$ are zero \cite{Carrell74}. On the other hand, he
also gave an example of a threefold $X$, a $\bP^1$-bundle over an abelian surface,
for which $c_1(X)^3\neq 0$. This suggests that one needs to treat varieties with
negative Kodaira dimension differently.

At the same time, Carrell's proof in the surface case starts by proving that a
surface admitting a nowhere vanishing holomorphic one form is necessarily minimal,
i.e., contains no $(-1)$-curves. Hence one might suggest the following.

\begin{proclaimspecial}{Wild Guess}\label{wg}
  If $X$ admits a nowhere vanishing holomorphic one form, then $X$ is minimal.
\end{proclaimspecial}

Unfortunately, we cannot expect this to hold in higher dimension: Let $X=A\times Y$
where $A$ is an abelian variety and $Y$ is arbitrary, or more generally let $X$ admit
a smooth morphism onto an abelian variety $A$ with general fibre $Y$. Then $X$ admits
nowhere vanishing holomorphic one forms, namely the ones pulled back from $A$, but if
$Y$ is not minimal, then neither is $X$.  The reason that this didn't happen for
surfaces is that every smooth curve is minimal.

So one may try the other part of the problem and ask whether the existence of a
nowhere vanishing holomorphic one form on a {\em minimal} variety $X$ implies that
$c_1(X)^{\dim X}=0$. For a minimal variety $X$, $K_X$ is nef, therefore $c_1(X)^{\dim
  X}\neq 0$ is equivalent to $K_X^{\dim X}>0$ which is equivalent to $X$ being of
general type.

We are also led to varieties of general type via a different path.  If $X$ admits a
nowhere vanishing holomorphic one-form, then \cite[Theorem 3.1]{Green-Laz87b} implies
that for generic $\sP\in \Pic ^0(X)$, one has $H^i(X,\Omega ^j_X\otimes \sP)=0$ for
all $i,j$. In particular, $\chi (X,\omega _X)=0$. %
On the other hand, when $X$ is of maximal Albanese dimension (i.e., $\dim X= \dim
\alb _X(X))$) and $\Alb(X)$ is simple, then X is a variety of general type if and
only if $\chi(X,\omega _X)>0$.

All of these considerations naturally lead to the following conjecture.

\begin{conj}\label{general}
  Let $X$ be a smooth projective variety of general type. Then $X$ does not admit a
  nowhere vanishing holomorphic one form.
\end{conj}


Finally observe, that once we restrict to varieties of general type, \myref{wg} does
not seem so wild anymore and one has a much more reasonable guess.

\begin{conj}[(Carrell)]\label{carrell}
  Let $X$ be a smooth projective variety of general type. If $X$ admits a nowhere
  vanishing holomorphic one form, then $X$ is minimal.
\end{conj}

\begin{rem}\label{lz}
  As mentioned above, this is known for surfaces and using the classification of
  extremal contractions one can easily see that it also holds for threefolds. This
  was explicitly checked in \cite[Lemma 2.1]{Luo-Zhang03}.
\end{rem}

Conjecture~\ref{general} has been confirmed for canonically polarized varieties
(i.e., whose canonical divisor is ample) in \cite{Zhang97} and for threefolds in
\cite{Luo-Zhang03}.  

An immediate consequence of this conjecture is that a variety of general type does
not admit any smooth morphisms onto an abelian variety.  For other applications the
reader is referred to \cite{Zhang97}.

In this article we first prove Conjecture~\ref{general} for smooth minimal varieties.

\begin{thm}[=Theorem~\ref{fornefandbig}]
\label{M1}
Let $X$ be a smooth projective minimal variety of general type. Then $X$ does not
admit a nowhere vanishing holomorphic one form.
\end{thm}

This completely confirms Conjecture~\ref{general} assuming Conjecture~\ref{carrell}.
Using \myref{lz} this also gives a new proof of the threefold case \cite[Theorem
1]{Luo-Zhang03}.

Using different methods than the ones used to prove \myref{M1}, we also confirm
Conjecture~\ref{general} for varieties whose Albanese variety is simple.

\begin{thm}[=Theorem~\ref{TS}]
\label{M2} 
Let $X$ be a smooth variety of general type. If its Albanese variety is simple, then
$X$ does not admit a nowhere vanishing holomorphic one form.
\end{thm}

\noin
\begin{demo}{Definitions and notation}
  Let $X$ be a proper variety. A line bundle $\sL$ on $X$ is called \emph{nef} if
  $\deg\left(\sL\resto C\right)\geq 0$ for every proper curve $C\subseteq X$.  $\sL$
  is called \emph{big} if the global sections of $\sL^m$ define a generically finite
  map for some $m>0$. $X$ is of \emph{general type} if $\omega_X$ is big.
  
  For $\theta\in\coh 0.X.\Omega_X.$, $Z(\theta)$ denotes the zero locus of $\theta$.
  
  Let $\sF$ be a torsion-free sheaf on $X$ and $\iota:U\into X$ the locus where $\sF$
  is locally free.  Then $\hat S^m(\sF)$ denotes the reflexive hull of the
  $m^{\text{th}}$ symmetric power of $\sF$, i.e., $\hat S^m(\sF)=\iota_*S^m(\sF\resto
  U)$.
\end{demo}

\Section{Smooth minimal models}
The main goal of this section is to prove the following.

\begin{proclaim}\label{main} Let $Y$ be a projective variety with only rational
singularities of dimension $n$, and let $\phi:\map X.Y.$ be a resolution of
singularities of $Y$.  Let $\foy=\im[\phi^*\Omega_Y\to \Omega_X]$.  Assume that \te s
a $\nattr\in\coh 0.X. \foy.$ such that the zero locus of $\nattr$ is empty.  Then for
any ample line bundle $\l$ on $Y$, $\coh n. Y. \l.=0$.
\end{proclaim}

Before we can prove this theorem we need some preparation.

Let $X$ be a smooth variety of dimension $n$.  Let $\Phi$ be the functor of regular
functions and $\Psi$ the functor of K\"ahler differentials, i.e., $\Phi_X=\ring X.$
and $\Psi_X=\Omega_{X}$. Then any $\nattr\in\coh 0.X.\Omega_{X}.$ induces a morphism
$\nattr_{X}:\map \Phi_X.\Psi_X.$. In fact it induces a morphism $\nattr_{X_i}:\map
\Phi_{X_i}.\Psi_{X_i}.$ via pull back for every $X_i$ that admits a morphism,
$\phi_i:X_i\to X$, to $X$. In other words $\nattr$ induces a natural transformation
from $\Phi$ to $\Psi$ in the category of $X$-schemes.  Then by \cite[{2.6,
  2.9}]{Kovacs04u} there exists a functorially defined $\Oxsui r.X.\in\ob(D(X))$ for
all $r\geq -1$ such that for every $p\in\N$ \te s a distinguished triangle,
\begin{equation}
\dist{\Oxsui p-1.X.}.{\Oxi p.}.{\Oxsui p.X.}.. \tag{$\star$}  
\end{equation}

Furthermore, $\Oxsui r.X.\simeq 0$ if $r> n-1$ and $\Oxsui
n-1.X.\simeq\ox$. 

Suppose $Z(\nattr)$ is empty. Then by \cite[{Appendix B.3.4}]{Fulton84} the Koszul
complex,
$$
0 \to \ring X. \stackrel{\land\nattr}{\longrightarrow} \Ox^1
\stackrel{\land\nattr}{\longrightarrow} \Ox^2
\stackrel{\land\nattr}{\longrightarrow} \dots
\stackrel{\land\nattr}{\longrightarrow} \Ox^{n-1}
\stackrel{\land\nattr}{\longrightarrow} \Ox^n \to 0,
$$ 

\noin induced by taking the wedge product with $\nattr$ is exact.  Let $\sE^{-1}=0$,
and $$\sE^i=\ker(\land\nattr):\map\Ox^{i}.\Ox^{i+1}.$$
for $i=1, \dots, n-1$. Then
$\Oxsui r.X.\simeq \sE^{r+1}$ for $r=1, \dots, n-2$. In
particular $\Oxsui 0.X.\simeq \ring X.$. 

\medskip

Next, results regarding the generalized De~Rham complexes are
summarized in the following theorem.

\begin{proclaim}{\rm (\cite{DuBois81}, \cite[{III.1.12, V.3.3,
      V.3.6, V.5.1)}]{GNPP88} }\label{derham} For every complex scheme $Y$ of
  dimension n there exists an $\Oyu^\cdot \in \ob(D_{filt}(Y))$ with the following
  properties.

\begin{roster}
\rostitem 
Let  $\phi_\cdot:\map X_\cdot.Y.$ be any hyperresolution of
$Y$. Then $\Oyu^\cdot\simeq\hypush\Oxh\cdot$.

\rostitem The definition is functorial, i.e., if $\phi :X\rightarrow Y$
is a morphism of complex schemes, then there exists a natural map
$\phi^{*}$ of filtered complexes $$ \phi^{*}:\Oyu^\cdot \rightarrow
R\phi_{*}\Oxu^{\cdot}.
$$
Furthermore, $\, \Oyu^\cdot \in \ob(D^{b}_{filt, coh}(Y))$ and if
$\phi$ is proper, then $\phi^{*}$ is a morphism in $D^{b}_{filt,
coh}(Y)$.\label{functorial}

\rostitem
Let\ $\, \Omega_{Y}^{\cdot}$ be the usual De~Rham complex of K\"ahler
differentials considered with the ``filtration b\^{e}te''. Then there
exists a natural map of filtered complexes 
$$
\Omega_{Y}^{\cdot}\rightarrow \Oyu^\cdot 
$$ 
and if $Y$ is smooth, it is a quasi-isomorphism.\label{usual}

\rostitem
Let $\Oyp = Gr^{p}_{F}\, \Oyu^\cdot [p]$. Then  $\Oyp
\simeq\hypush\Oxph$ for any hyperresolution $\phi_\cdot:\map
X_\cdot.Y.$. \label{independent}

\rostitem
If $Y$ is projective and $\l$ is an ample line bundle on $Y$, then $$
{\mathbb H}^{q}(Y,\Oyp\otimes\l)=0 \qquad \text{for } p+q > n.
$$\label{hypkan}

\end{roster}
\end{proclaim}

\noin
To extend the definition of $\Oxsui p.X.$ to singular varieties we
need the following.

\begin{lem}
  Let $\phi_\cdot:\map X_\cdot.Y.$ be a hyperresolution of\ $Y$. Let
  $\foy$ be defined as in \myref{main}. Let $\nattr\in\coh 0.X_0.
  \foy.$ and $\nattr_{X_i}:\ring X_i.\to\Omega_{X_i}$ the morphism
  induced by the section $\nattr$.  Then $\rpforward{\phi_\cdot}\Oxsui
  r.{X_\cdot}.$ is independent of the hyperresolution chosen.
\end{lem}

\proof
Let $\alpha$ be a morphism of hyperresolutions. 
$$
\CD
X'_\cdot @>\alpha>> X''_\cdot \\
@V\varepsilon'_\cdot VV @VV\varepsilon''_\cdot V \\
X @>>\id_X > X 
\endCD
$$

\noin Then by \cite[{4.1}]{Kovacs04u} \te s a commutative diagram:
$$
\CD
R\varepsilon'_\cdot\Oxsui p-1.X'_\cdot. @>>>
R\varepsilon'_\cdot\Omega^p_{X'_\cdot} @>>>  
R\varepsilon'_\cdot\Oxsui p.{X'_\cdot}. @> +1 >>  \\
@VVV @VVV @VVV \\
R\varepsilon''_\cdot\Oxsui p-1.X''_\cdot. @>>>
R\varepsilon''_\cdot\Omega^p_{X''_\cdot} @>>>  
R\varepsilon''_\cdot\Oxsui p.{X''_\cdot}. @> +1 >>  \\
\endCD
$$

\noin Now $R\varepsilon'_\cdot\Omega^p_{X'_\cdot}\simeq\Oyp\simeq
R\varepsilon''_\cdot\Omega^p_{X''_\cdot}$ by \myref{independent}, and the statement
follows from \cite[2.1.4]{DuBois81} and ($\star$) by descending induction on $p$.
\endproof

\begin{defn}
Let $Y$ be a variety of dimension $n$. Let $\phi_\cdot:\map
X_\cdot.Y.$ be a hyperresolution of $Y$ and let $\nattr\in\coh
0.X_0. \foy.$.  We define $\Oxsui r.Y.=\rpforward{\phi_\cdot}\Oxsui
r.{X_\cdot}.$ for $r\geq -1$. By the lemma, this is independent of the
hyper\-resolution chosen, in particular if $Y$ is smooth, it agrees
with the previous definition of $\Oxsui r.Y.$.
\end{defn}

\medskip
\noin
{\it Proof of Theorem~\ref{main}.}
By ($\star$) \te s a distinguished triangle,
\begin{align}
  \dist{\Oxsui p-1.Y.}.{\Oyp}.{\Oxsui p.Y.}..\tag{$\star\star$}
\end{align}

\noin for every $p\in\N$, so by \myref{hypkan}, $$\hypcoh n-p.Y.{\Oxsui
  p.Y.\tensor\l}.\to \hypcoh n-(p-1).Y.{\Oxsui p-1.Y.\tensor\l}.$$
is surjective for
all $p$, and then $$\hypcoh 0.Y.{\Oxsui n.Y.\tensor\l}.\to\dots\to \hypcoh
n.Y.{\Oxsui 0.Y.\tensor\l}.$$
is also surjective. Now $\hypcoh 0.Y.{\Oxsui
  n.Y.\tensor\l}.=0$ since $\Oxsui n.Y.=0$, so we obtain that $$\hypcoh n.Y.{\Oxsui
  0.Y.\tensor\l}.=0.$$

On the other hand, the previous observation in the case $Z(\nattr)=\emptyset$,
\myref{functorial}, \myref{usual}, ($\star\star$) and \cite[4.1]{Kovacs04u} implies
that the following diagram is commutative:
$$
\CD
\ring Y. @>>> \underline{\Omega}^0_Y @>>> \Oxsui 0.Y. \\
@VV\rho V @VVV @VVV \\
\rpforward{\phi}\ring X. @>\simeq >>
\rpforward{\phi}\underline{\Omega}^0_{X} @>\simeq >> 
\rpforward{\phi}\Oxsui 0.X..
\endCD
$$
Now $\rho$ has a left inverse, and hence in turn the morphism
$\ring Y.\to \Oxsui 0.Y.$ has a left inverse. Finally that implies
that $\coh n. Y.  \l.\to\hypcoh n.Y.{\Oxsui 0.Y.\tensor\l}.=0$ is
injective.  \qed

\begin{proclaim}\label{fornefandbig}
Let $X$ be a smooth projective variety such that $\omega_X$ is nef and
big, and let $\nattr\in H^0(X,\Omega_X)$. Then $Z(\nattr)\neq\emptyset$.
\end{proclaim}

\proof
Consider the Albanese morphism, $h:X\to A$ and let
$$g:\map{Y=\proj A.h_*\omega_X^m.}.A.$$
and $\phi:X\to Y$ the induced natural
morphism. Note that by construction $h$ factors as $g\circ \phi$. 

$Y$ has rational singularities by \cite{Elkik81} (cf.\ \cite{Kovacs00b}), and
$\omega_Y$ is a line bundle by the Basepoint-free theorem \cite[Theorem~3.3]{KM98}
(cf.\ \cite{Reid83}). In particular, $Y$ is Gorenstein.

Next let $\eta\in\coh 0.A.\Omega_A.$ \st\ $\nattr=h^*\eta$. It follows that
$\nattr\in\coh 0.X. \foy.$, and hence by Theorem~\ref{main}, if $\oy$ is ample, then
$Z(\nattr)\neq\emptyset$.  Therefore it is enough to prove that $\oy$ is ample.

Consider the canonical model of $X$: $$Z=\proj._{\text{\tiny $m\geq
    0$}} H^0(X,\omega_X^m)..$$
We want to prove that $Y\simeq Z$.

By construction there is a natural morphism $\psi: X\to Z$ that
factors through $\phi$. Let $\eta:Y\to Z$ be the induced morphism such
that $\psi=\eta\circ\phi$.

For all $a\in A$, $\omega_{Y_a}$ is ample. On the other hand, for any
curve $C$ contained in a fibre of $\eta$, $\omega_Y=\eta^*\omega_Z$ is
trivial on $C$. Hence $Y_a$ intersects every fiber of $\eta$ in a
zero dimensional subscheme.
 
Let $E$ be a component of the exceptional locus of $\eta$.  By
\cite[Theorem~2]{Kawamata91}, $E$ is covered by rational curves that are contracted
by $\eta$. By the above observation these rational curves cannot be contained in any
of the $Y_a$.  Since $A$ does not contain any rational curves we conclude that $E$
must be empty.

Thus $Y\simeq Z$, in particular $\omega_Y$ is ample.
\endproof


\Section{Varieties whose Albanese variety is simple}
\begin{thm} \label{TS} Let
  $X$ be a smooth variety of general type and $\alb_X:X\to A\leteq\Alb(X)$ its
  Albanese morphism.  If $A$ is simple, then any holomorphic one-form $\theta\in
  H^0(X,\Omega ^1 _X)$ has a non-empty zero set.
\end{thm}

We are going to study the Albanese morphism of $X$ and employ different strategies
depending on whether it is surjective or not.

\sh{\sc Case I: $\alb_X$ is not surjective}
\begin{prop} 
  Let $Z\subsetneq A$ be a proper closed subvariety of the abelian variety $A$.  If
  $A$ is simple, then for every holomorphic one-form $\theta\in H^0(A,\Omega ^1 _A)$,
  $\theta\resto Z$ has a non-empty zero set.
\end{prop}

\begin{proof} 
  Let $\cW \subset H^0(A,\Omega ^1 _A)$ be the set of those holomorphic one-forms
  $\theta\in H^0(A, \Omega ^1 _A)$ such that $\theta\resto Z$ vanishes at some point
  $z\in Z$.  It is easy to see that $\cW$ is closed and so it suffices to show that
  $\cW$ is dense in $H^0(A, \Omega ^1 _A)$.
  
  Let $r=\dim Z$ and $Z_0$ the set of smooth points of $Z$. For any $z\in Z_0$, one
  has that the tangent space $T _z(Z)^\perp\cong \bC ^{g-r}\subset T _z(A)^\vee \cong
  H^0(A,\Omega ^1 _A)\cong \bC ^{g}$.  Let $\cZ \subset \cP ^{g-1}$ be the closure of
  the image $\cZ _0$ of the corresponding projective bundle $\cP :=\PP (T
  (Z_0)^\perp)$ under this map. One sees that if $\cZ =\PP ^{g-1}$, then $\cW$ is
  dense in $H^0(A, \Omega ^1 _A)$.
  
  Suppose that $\cZ\neq\PP^{g-1}$, i.e., $\dim\cZ<g-1$.  Let $p\in \cZ$ be a general
  point, then $\dim \cZ < g-1$ implies that the corresponding fiber $\cP _p$ is
  positive dimensional.  Consider now the projection $\pi:\cP \to Z_0$ and the
  subvariety $Z _p$ given by the closure of $\pi (\cP _p) \cap Z_0\subset A$, one has
  $$\dim Z _p=g-1-\dim \cZ>0.$$
  For general $x\in Z_p$, one has for $L_p$ the line
  corresponding to $p$ that $L_p\subset T_x(Z)^\perp $ and so $T_x(Z_p)\subset
  T_x(Z)\subset H_p:=L_p^\perp$.  It follows that $Z_p$ generates a proper abelian
  subvariety $A_p\subsetneq A$. This is a contradiction, so $\cW=H^0(A,\Omega ^1
  _A)$.
\end{proof}

\begin{cor} \label{P2} Let $X$ be a projective variety and
  $\alpha:X\to A$ a morphism to a simple abelian variety.  If $Z:=\alpha(X)\ne A$,
  then every holomorphic one-form
  $$\theta\in \alpha^*H^0(A,\Omega ^1 _A)\subset H^0(X, \Omega^1_X)$$
  has a non-empty
  zero set.
\end{cor}

\subheading{\sc Case II: $\alb_X$ is surjective}\hfil\penalty -10000\vspace{-12pt}
\begin{notation}\label{setup}
  Let $X$ be a projective variety and $\alpha:X\to A$ a surjective morphism to an
  abelian variety.  Let $\Delta\subset A$ be the locus where $\alpha$ is not smooth.
\end{notation}

\begin{prop} \label{P3} 
  Under the assumptions of \myref{setup} assume that $\Delta$ contains an irreducible
  divisor $D$ of general type. Then every holomorphic one-form $\theta\in
  \alpha^*H^0(A, \Omega ^1 _A)\subset H^0(X, \Omega ^1 _X)$ has a non-empty zero set.
\end{prop}
\begin{proof}
  Consider $\cW \subset H^0(A,\Omega ^1 _A)$ the set of those holomorphic one-forms
  $\theta\in H^0(A,\Omega ^1 _A)$ such that $\alpha^*\theta\in H^0(X,\Omega ^1 _X)$
  vanishes at some point $x\in X$.  As above, $\cW$ is closed and so it suffices to
  show that it is dense.
  
  Consider $D^0\subset D$ a (non-empty) open set such that for all $z\in D^0$ there
  is a point $x\in X_z$ with $\rank (d\alpha_x )= g-1$ (cf. \cite[III.10.6]{Ha77})
  and $D$ is smooth at $z$. Let $x_1,...,x_{g}$ be local coordinates of $A$ at $z$
  such that $D$ is defined by $x_g=0$ and $\theta=\theta _z \in H^0(A,\Omega ^1_A)$
  such that $\theta(z)=dx_g$.  Then, $\theta$ spans the subspace $T_z(D)^\perp\subset
  H^0(A,\Omega ^1 _A)$ and $\theta\resto D$ vanishes at $z$ and $\alpha^*\theta$
  vanishes at some point $x\in X$ such that $a(x)=z$ (in fact at any such point with
  $\rank (d\alpha_x )= g-1$).  Since $D$ is of general type, by \cite{Gri-Har79}
  (cf.\ \cite[(3.9)]{Mori87}), its Gauss map is generically finite and so one sees
  that the set $\{\theta_z|z\in D^0 \}\subset \cW$ is dense in $H^0(A,\Omega ^1 _A)$.
  (Reasoning as in the previous proposition, we have that $\cP \cong Z_0$ and $\cP
  \to \PP ^{g-1}$ is generically finite and so it is dominant.)
\end{proof}

\begin{prop}\label{P4}
  Under the assumptions of \myref{setup} assume that there exists a positive integer
  $m$ such that $\alpha_*(\omega _{X/A}^{m})$ is big. Then $\alpha$ is not smooth
  in codimension one, i.e., $\Delta$ contains a divisor.
\end{prop}

\begin{proof} 
  Since $\alpha_*(\omega _{X/A}^{m})$ is big, for any ample line bundle $\sH$ on $A$
  there exists an integer $a>0$ such that $\hat{S}^a(\alpha_*(\omega _{X/A}^{\ot
    m}))\ot \sH^{-1}$ is big. Let $m_k:A_k\simeq A\to A$ be multiplication by an
  integer $k$, so $m_k$ is an \'etale map such that $m_k^*\sH=(\sH_k)^k$ with $\sH_k$
  an ample line bundle on $A_k$.  Let $g=\dim A$, $r=\dim X-\dim A$, and
  $k=3r(g-1)ma$. Further let $H_k$ be a divisor on $A_k$ such that $\sO_A(H_k)\simeq
  \sH_k$ and finally let
  $$\alpha':X':=X\times _A A_k\to A_k=A.$$
  Then
  $$m_k^*(\hat{S}^a(\alpha_*(\omega _{X/A}^{m}))\ot \sH^{-1})=
  \hat{S}^a({\alpha'}_*(\omega _{X'/A}^{m}))\ot ( {\sH_k}^{3r(g-1)})^{-ma}$$
  is big
  (cf.  \cite[(5.1.1) (d)]{Mori87}) and hence ${\alpha'}_*(\omega _{X'/A}^{m})\ot (
  {\sH_k}^{3 r(g-1)})^{-m}$ is also big.  Since $H_k$ is ample, $3H_k$ is very ample.
  Let $C$ be a curve obtained by intersecting $g-1$ general elements in $\vert 3
  H_k\vert$ and $\sA :=\sH_k^{3r(g-1)}\resto C$.  Then
  $$\deg \omega_C=(g-1)(3H_k)^g \quad\text{and}\quad \deg \sA
  =3r(g-1)H_k\cdot(3H_k)^{(g-1)}= r\deg \omega_C.$$
  Let $Y=(\alpha')^{-1}(C)$ and
  $h=\alpha'\resto Y:Y\to C$.  If $\alpha$ is smooth in codimension one, then
  $\alpha'$ is smooth in codimension one and so $h$ is smooth.  Since
  $$\leftresto\left( {\alpha'}_*(\omega _{X'/A}^{m})\ot ( {\sH_k}^{3r(g-1)})^{-m}
  \right)\rightresto C= h_*(\omega _{Y/C}^{m})\ot \sA ^{-m},$$
  it follows that
  $h_*(\omega _{Y/C}^{m})\ot \sA ^{-m}$ is also big and hence ample.  By
  \cite[Proposition 4.1]{Vie-Zuo01} (with $\delta =0$, $s=0$), one has
  $$\deg\sA <\dim (Y/C)\deg \omega _C=r \deg \omega _C.$$
  This is impossible and
  hence $\alpha$ is not smooth in codimension one.
\end{proof}

\begin{proof}({of Theorem \ref{TS}})
  Since $A$ is simple, by \cite[10.9]{Ueno75} (cf. \cite[Theorem 3.7]{Mori87}) any
  proper subvariety of $A$ is of general type. By \myref{P2} we may assume that
  $\alb_X:X\to A$ is surjective. Then by \myref{P3}, we may assume that $\alb_X$ is
  smooth in codimension one (again using \cite[10.9]{Ueno75} to see that every
  divisor is of general type).
  
  Now let $X\to Z\to A$ be the Stein factorization of $\alb_X$. Then $Z\to A$ is
  smooth and hence \'etale in codimension one, so $Z$ is birational to an abelian
  variety. It follows that $Z$ is birational to $A$ and $\alb_X:X\to A$ is an
  algebraic fiber space.  Since $X$ is of general type, $(\alb_X)_*(\omega
  _{X/A}^m)=(\alb_X)_*(\omega _X^m)$ is big for some $m>0$, but by \myref{P4} this is
  impossible.
\end{proof}

\bibliographystyle{$HOME/tex/TeX_input/myalpha} 
\bibliography{$HOME/tex/TeX_input/ref} 

\end{document}